\documentclass[letterpaper, 10 pt, conference]{ieeeconf}  % Comment this line out if you need a4paper

\IEEEoverridecommandlockouts                              % This command is only needed if 
\usepackage{amsmath,amssymb,amsfonts}
\usepackage{graphicx}
\usepackage{textcomp}
\def\BibTeX{{\rm B\kern-.05em{\sc i\kern-.025em b}\kern-.08em
    T\kern-.1667em\lower.7ex\hbox{E}\kern-.125emX}}

\usepackage{bm}
\usepackage{algorithm}
\usepackage{algpseudocode}
\usepackage{cases}
\usepackage{amsmath,mathtools}
\usepackage{nccmath}
\usepackage{etoolbox}
\preto\subequations{\ifhmode\unskip\fi}

\makeatletter
\newenvironment{paradigm}[1][htb]{%
    \renewcommand{\ALG@name}{Paradigm}% Update algorithm name
   \begin{algorithm}[#1]%
  }{\end{algorithm}}
\makeatother

\title{\LARGE \bf
Privacy-Preserving Decentralized Multi-Agent Cooperative Optimization -- Paradigm Design and Privacy Analysis 
}

\author{Xiang Huo$^{1}$ and Mingxi Liu$^{1}$% <-this % stops a space

\thanks{$^{1}$X. Huo and M. Liu are with the Department of Electrical and Computer Engineering,
University of Utah, 50 S Central Campus Drive, Salt Lake City, UT, 84112, USA
{\tt{\small\{xiang.huo, mingxi.liu\}@utah.edu}}}%
}

\begin{document}

\maketitle
\thispagestyle{empty}
\pagestyle{empty}

%%%%%%%%%%%%%%%%%%%%%%%%%%%%%%%%%%%%%%%%%%%%%%%%%%%%%%%%%%%%%%%%%%%%%%%%%%%%%%%%
\begin{abstract}
Large-scale multi-agent cooperative control problems have materially enjoyed the scalability, adaptivity, and flexibility of decentralized optimization. However, due to the mandatory iterative communications between the agents and the system operator, the decentralized architecture is vulnerable to malicious attacks and privacy breach. Current research on addressing privacy preservation of both agents and the system operator in cooperative decentralized optimization with strongly coupled objective functions and constraints is still primitive. To fill in the gaps, this paper proposes a novel privacy-preserving decentralized optimization paradigm based on Paillier cryptosystem. The proposed paradigm achieves ideal correctness and security, as well as resists attacks from a range of adversaries. The efficacy and efficiency of the proposed approach are verified via numerical simulations and a real-world physical platform.

%Currently, there are few studies addressing decentralized algorithms from the aspect of strongly coupled objective functions and constraints, and jointly with privacy preservation of both the agents and the system operator. 

\end{abstract}

%%%%%%%%%%%%%%%%%%%%%%%%%%%%%%%%%%%%%%%%%%%%%%%%%%%%%%%%%%%%%%%%%%%%%%%%%%%%%%%%
\section{Introduction}

Decentralized optimization has attracted remarkable attentions in large-scale cooperative multi-agent control owing to its scalability, adaptivity, and flexibility, and has been applied to a variety of applications ranging from multi-user optimization \cite{koshal2011multiuser}, electric vehicle charging control \cite{nimalsiri2019survey}, evolutionary computation \cite{dorronsoro2011improving}, to electricity market design \cite{morstyn2018designing}. However, in any decentralized cooperative optimization architecture, the mandatory iterative information exchange between the agents and the system operator (SO) would pose potential privacy risks, including the leakage of intermediate decision variables and the exposure of private coefficients of both the agents and the SO \cite{lu2018privacy}. 

Differential privacy (DP) has been commonly used to prevent privacy leakage \cite{han2016differentially, zhang2016dynamic}. DP achieves privacy preservation through masking the exchanged sensitive information by adding deliberate noises \cite{han2016differentially}. However, DP-based strategies inevitably suffer from accuracy compromise caused by the introduced noises. In \cite{zhang2016dynamic}, an optimal dynamic DP mechanism with guidelines to choosing privacy parameters was proposed to minimize the accuracy deterioration. Nonetheless,  a fundamental  trade-off between privacy and accuracy universally exists in DP-based strategies \cite{zhang2018enabling}. In resolving this issue, Wang \cite{wang2019privacy} proposed an approach in which each agent decomposes its state into two sub-states and only reveals one to its neighbors. This noise-free decomposition technique enables privacy preservation without compromising accuracy, however can only deal with decoupled optimization problems.  

Another approach to achieving privacy preservation with ideal accuracy is through cryptography. Cryptography-based strategies can be classified into two types: symmetric and asymmetric \cite{johnson2019security}. In \emph{symmetric} cryptography, which also referred to as private key cryptography, each pair of partners share the same private key and keep it secret from others. The decryption and encryption will adopt the same private key. Standard private key cryptography algorithms include Data Encryption Standard (DES), Advanced Encryption Standard (AES), SingleMod, etc. Comparing with asymmetric cryptography, symmetric cryptography requires less computational resource, while lacking reliability and security due to the key selection, distribution, and storage. In \emph{asymmetric} cryptography, which also referred to as public key cryptography, each  agent has two pairs of key: a public key and a private key. The public key is available to anyone and the private key is kept secret. Anyone can use the public key to encrypt a plaintext to a ciphertext, but only the private key host could decrypt the ciphertext. Asymmetric cryptography has an enhanced security compared to symmetric cryptography and is scalable to the large-scale agent population size and frequent information exchange. Examples of public key cryptography include Rivest–Shamir–Adleman (RSA), Digital Signature Algorithm (DSA), Paillier cryptosystem, etc. This paper aims to integrate public key cryptography into decentralized optimization algorithm design.

Due to the iterative nature of decentralized optimization, incorporating cryptography into the algorithm design mandates arithmetic operations over ciphertext. Homomorphic cryptosystems allow arithmetic operations over ciphertexts, and when decrypting the outcome, the decrypted value matches the result of operations performed on the plaintexts. Partially homomorphic cryptosystem, e.g., RSA, can execute addition or multiplication operations, while fully homomorphic cryptosystem, e.g., SingleMod, can execute both addition and multiplication operations. In \cite{hadjicostis2020privacy}, a ratio consensus algorithm was designed based on partially homomorphic encryption. However, it is only applicable to integers and obliged with an assumption of the presence of a trusted agent. Lu \emph{et al.} \cite{lu2018privacy} proposed a private key based homomorphic encryption that is applicable to generic distributed projected gradient-based algorithms. However, the private key based algorithm is not semantically secure. This paper is set to develop a semantically secure paradigm with additional assumptions eliminated, e.g., without the existence of a trusted agent.

In this paper, we work towards designing a privacy-preserving decentralized optimization paradigm for strongly coupled cooperative optimization problems. The contribution of this paper is three-fold: (1) We propose a novel decentralized privacy-preserving paradigm that is semantically secure towards various types of adversaries. Compared with DP-based strategies, we allow the participants to directly exchange their private information without any perturbation and achieve the exact optimality; (2) The privacy preservation of all participants (agents and the SO) is guaranteed and rigorously proved for a class of strongly coupled optimization problems. Moreover, the proposed paradigm naturally accomplishes the coefficient assignment task and achieves enhanced security towards the convergence of the decision variables; (3) We build a real-world physical platform with Raspberry Pi boards to demonstrate the industrial values of the proposed privacy-preserving paradigm.

\section{Problem Formulation}

Consider a coupled decentralized cooperative optimization problem involving $n$ agents, where the $i$th agent holds a local cost function $f_i(x_i)$ with a local decision variable $x_i$ that is locally constrained by a convex set $\mathbb{X}_i$. In the presence of a global objective function and a global constraint, the cooperative optimization problem can be formulated as 
\begin{subequations} \label{4}
\begin{align} 
& \underset{x}{\text{min}} & & \frac{1}{2} \|\sum_{i=1}^{n} A_{ui} x_i + c \|^2_2
    +\sum_{i=1}^{n}f_i(x_i) \label{obj} \\
& \: \text{s.t.} & & x_{i} \in \mathbb{X}_{i}, \quad \forall i=1,\cdots,n, \\
& & & \sum_{i=1}^{n} A_{gi} x_i + d \leq 0 , \label{global_con}
\end{align}
\end{subequations} 
where the quadratic term in \eqref{obj} is the coupled global objective, \eqref{global_con} denotes the coupled global constraint,  $A_{ui}$ and $A_{gi}$ are local coefficients of the network that are assigned to the $i$th agent, and $c$ and $d$ are aggregated global vectors containing the network information that are assigned to the SO. Note that the SO also has access to $A_{ui},A_{gi}, \forall i\in \mathcal{N}$, where $\mathcal{N}$ denotes the set of agents. For example, for the electric vehicle charging control problem in \cite{liu2017decentralized}, $c$ represents the baseline load profiles and $d$ represents the voltage bounds; for the resource allocation problem in \cite{xu2017distributed},
$d$ represents the resource capacity limits.

Eqn. \eqref{4} is strongly coupled through the quadratic term in the objective function and the globally coupled constraints. Primal-dual methods e.g., regularized primal-dual subgradient (RPDS) \cite{koshal2011multiuser} and shrunken primal-dual subgradient (SPDS) \cite{liu2017decentralized}, can be used to solve this class of problem in a decentralized way. To use any primal-dual method, one firstly calculate the relaxed Lagrangian of \eqref{4} as
\begin{equation} \label{5}
    \mathcal{L}(x,\lambda)= F(x)+\lambda^{\mathsf{T}}l(x), 
\end{equation}
where $x=[x_1^{\mathsf{T}} \cdots x_n^{\mathsf{T}}]^{\mathsf{T}}$, $F(x) \triangleq \frac{1}{2} \|\sum_{i=1}^{n} A_{ui} x_i + c \|^2_2
    +\sum_{i=1}^{n}f_i(x_i)$, $l(x) \triangleq \sum_{i=1}^{n} A_{gi} x_i + d$, and $\lambda$ denotes the dual variable. Consequently, the subgradients of $\mathcal{L}(x,\lambda)$ \emph{w.r.t.} $x_i$ and $\lambda$ can be calculated as
\begin{subequations}
\begin{align}
    \nabla_{x_i} \mathcal{L}(x,\lambda) &{=} {A_{ui}^{\mathsf{T}}} (\sum_{i=1}^{n} A_{ui} x_i {+} c )
    {+}\nabla_{x_i} f_i(x_i) {+}  {A_{gi}}^\mathsf{T}\lambda, \label{4as}\\
     \nabla_{\lambda} \mathcal{L}(x,\lambda) & {=} \sum_{i=1}^{n} A_{gi} x_i + d.    \label{4bs}
\end{align}
\end{subequations}
Note that both RPDS and SPDS can be integrated into the privacy-preserving paradigm proposed in Section \ref{Section_Paradigm}. We take SPDS, which has higher accuracy and convergence rate, for illustration of the paradigm design. The primal and dual updates in SPDS follow
\begin{subequations} \label{SPDS_update}
\begin{align}
x_{i}^{k+1} &=\Pi_{\mathbb{X}_{i}}\left( \frac{1}{\tau_x} 
\Pi_{\mathbb{X}_{i}}\left( \tau_x x_{i}^{k}-\alpha_{i} \nabla_{x_{i}} \mathcal{L}\left(x^{k}, \lambda^{k}\right)\right)\right), \label{11a}\\
\lambda^{k+1} &= \Pi_{\mathbb{D}}\left( \frac{1}{\tau_{\lambda}}\Pi_{\mathbb{D}}\left(
\tau_{\lambda}\lambda^{k}+\beta \nabla_{\lambda} \mathcal{L}\left(x^{k}, \lambda^{k}\right)\right)\right),\label{11b}\
\end{align}
\end{subequations}
where $k$ is the iteration index, $\tau_x$ and $\tau_{\lambda}$ are the shrunken parameters for the primal and dual updates, respectively, $\alpha_i$ and $\beta$ are the primal and dual step sizes, respectively, and $\mathbb{D}$ is a convex constraint set for the dual variable. The details of SPDS can be referred to \cite{liu2017decentralized}. It can be seen from  \eqref{SPDS_update} that both the primal and dual updates require the decision variables from all agents. Thus, intermediate decision variables of the agents must be transmitted to the SO, putting agents' privacy at risk. To resolve this issue, we aim to develop a decentralized 
privacy-preserving paradigm that protects the privacy of both the agents and the SO in solving problem \eqref{4}.

\section{Paradigm Design} \label{Section_Paradigm}

\subsection{Paillier Cryptosystem}

Paillier cryptosystem is  based on public key cryptography, and it takes the advantage of additively homomorphic, semantically secure encryption, and efficient decryption. Due to the page limit, we only present the properties that are necessary for establishing our design architecture. In Paillier cryptosystem, the participant generates two sets of keys, i.e., the public key $(p,g)$ and the private key $(\eta,\mu)$, then it publicizes the public key while keeping the private key to itself. Anyone with the public key can encrypt an integer plaintext $m_z$ by
\begin{equation}
c=g^{m_z} r^{p} \bmod p^{2},
\label{1}
\end{equation}
where $r$ is a random integer, $c$ denotes the ciphertext, and $\bmod$ denotes the modular operation. The ciphertext can be decrypted with the private key by
\begin{equation}
 \tilde{m}_z =  \left\lfloor\frac{\left(c^{\eta} \bmod p^{2}\right)-1}{p}\right\rfloor \mu \bmod p,
 \label{2}
\end{equation}
where $\tilde{m}_z$ denotes the decrypted message and $\lfloor \cdot \rfloor$ denotes the floor of a real number. Note that the coefficients and variables (e.g., $A_{ui}$ and $x_i$) in the decentralized calculations are real numbers, however $m_z$ has to be an integer in any cryptosystem.
Therefore, a real number $r_e$ is transformed into an integer $m_z$ by $m_z = 10^{\sigma} r_e$, where $\sigma$ denotes the preserved decimal fraction digits. A ciphertext can be transformed back to a real number by \cite{lu2018privacy}
\begin{align}
T_{\sigma, p}(\tilde{m}_z)= \begin{cases}
     \tilde{m}_z/ 10^{\sigma}, &\text { if } 0 \leq \tilde{m}_z \leq(p-1) / 2, \\ 
    (\tilde{m}_z-p) / 10^{\sigma}, &\text { if }(p+1) / 2 \leq \tilde{m}_z<p.
    \end{cases}
    \label{12}
\end{align}
The Paillier cryptosystem is additively homomorphic and multiplicatively semi-homomorphic, indicating that
\begin{subequations} \label{fully_homorphic}
\begin{align}
\mathcal{D}(\prod_{\ell=1}^{n} \mathcal{E}(m_{z\ell},p,g,r_{\ell})) &=\sum_{\ell=1}^{n} m_{z\ell}, \label{4a} \\
\mathcal{D}(\mathcal{E}(m_{z1},p,g,r_{1})^{m_{z2}})&=m_{z1} m_{z2}, \label{4b} 
\end{align}
\end{subequations}
where $z_{r\ell}$ and $r_{\ell}$ are the $\ell$th plaintext message and $\ell$th random integer, respectively, and $\mathcal{E}(\cdot)$ and $\mathcal{D}(\cdot)$ denote the encryption and decryption processes, respectively. Note that in \eqref{4b}, the multiplicatively semi-homomorphic property comes from the need of a plaintext power factor $m_{z2}$ for the purpose of getting the multiplication of two plaintexts, i.e., $m_{z1}m_{z2}$.
In the rest of the paper, we regard the real numbers have been transformed into integers for simplicity.

\subsection{Privacy-Preserving Paradigm Design}

\begin{figure}[!htbp]
% \vspace*{-2mm}
    \centering
 \includegraphics[width=0.5\textwidth,trim = 0mm 0mm 0mm 0mm, clip]{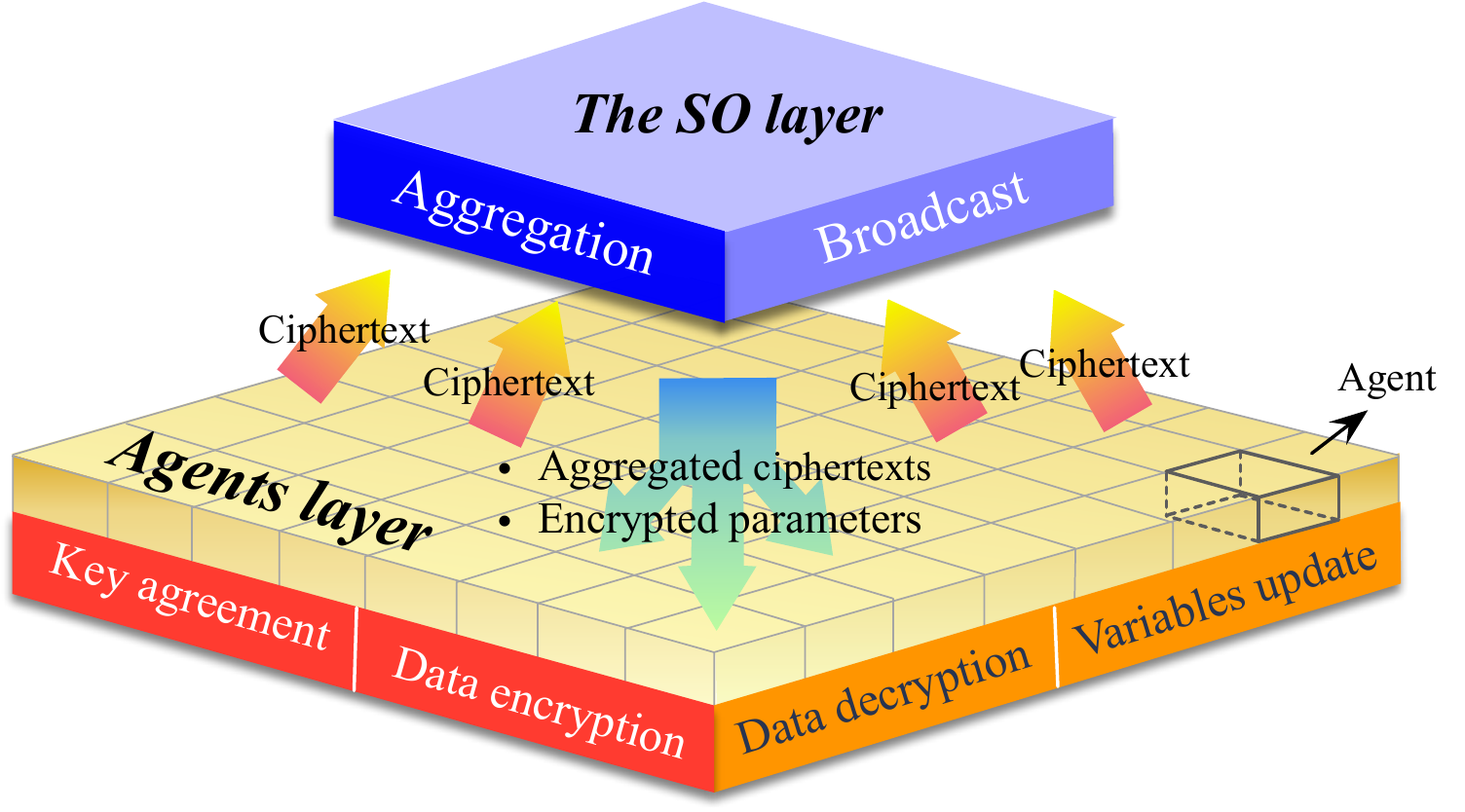}
    \caption{Decision variable and dual variable updating process.}
    \label{update}
\end{figure}

In this section, we propose a novel cryptography-based privacy-preserving decentralized optimization paradigm for problem \eqref{4}. Fig. \ref{update}  presents the information exchange and updating architecture of the paradigm. Before the iteration begins, the agents agree on a set of public key $(p,g)$ and private key $(\eta,\mu)$, and broadcast the public key to the SO while holding the private key to themselves. During each iteration of any primal-dual based decentralized algorithm,
the SO randomly generates two sets of parameters, satisfying 
\begin{subequations}
\begin{align}
    \gamma^k_1+\gamma^k_2+\cdots+\gamma^k_n &= 1, \label{8a}\\
    \upsilon^k_1+\upsilon^k_2+\cdots+\upsilon^k_n &= 1,
\end{align}
\end{subequations}
and sends $\mathcal{E}(\gamma_i^kc)$ and $\mathcal{E}(\upsilon_i^kd)$ to the $i$th agent. The $i$th agent encrypts $ A_{ui} x_i^k + \gamma_i^kc$ and $ A_{gi} x_i^k + \upsilon_i^kd$ using the public key. Then, the agents send $\mathcal{E}(A_{ui} x_i^k + \gamma_i^kc)$ and $\mathcal{E}(A_{gi} x_i^k + \upsilon_i^kd)$, $\forall i \in \mathcal{N}$ to the SO. Therefore, only encrypted messages are exchanged between the agents and the SO. The SO then calculates the multiplication of the received ciphertexts based on the addtively homomorphic property of the Paillier cryptosystem and broadcasts 
$\mathcal{E}(\sum_{i=1}^{n} A_{ui} x_i^k + c)$ and $\mathcal{E}(\sum_{i=1}^{n} A_{gi} x_i^k + d)$ to all the agents. The agents receive and decrypt $\mathcal{E}(\sum_{i=1}^{n} A_{ui} x_i^k + c)$ and $\mathcal{E}(\sum_{i=1}^{n} A_{gi} x_i^k + d)$ with the private key $(\eta,\mu)$, then convert the integers to real numbers. Finally, the $i$th agent updates the primal variable $x_i$ by \eqref{11a} and the dual variable $\lambda$ by \eqref{11b} over plaintext. The detailed steps of the proposed decentralized  privacy-preserving paradigm is presented via Paradigm \ref{paradigm_1}, and the rigorous privacy analysis is provided in Section \ref{Section_analysis}.

\begin{paradigm} 
\caption{Cryptography-based decentralized privacy-preserving multi-agent cooperative optimization}
\begin{algorithmic}[1]
% \Procedure{CH\textendash Election}{}

\State All agents agree on the public key $(p,g)$ and private key $(\eta,\mu)$. Then, agents broadcast $(p,g)$ to the SO and keep $(\eta,\mu)$ private.

\State All agents initialize primal and dual variables, tolerance $\epsilon_0$, iteration counter $k{=}0$, and maximum iteration $k_{max}$.

\While{ $\epsilon > \epsilon_0$ and $k < k_{max}$}

\State SO generates a set of parameters satisfying $\sum_{i=1}^n \gamma_i^k = 1$ and 
$\sum_{i=1}^n \upsilon_i^k {=} 1$, then uses the public key $(p,g)$ for encryption and sends $\mathcal{E}(\gamma_i^kc)$ and $\mathcal{E}(\upsilon_i^kd)$  to the $i$th agent.

\State  All agents encrypt $ A_{ui} x_i^k + \gamma_i^kc$ and $ A_{gi} x_i^k + \upsilon_i^kd$ using Paillier encryption in \eqref{1} with public key $(p,g)$, then send $\mathcal{E}(A_{ui} x_i^k + \gamma_i^kc)$ and $\mathcal{E}(A_{gi} x_i^k + \upsilon_i^kd)$ for $i=1,\ldots,n$ to the SO.

\State  SO firstly collects the ciphertexts from all agents, then multiplies the received ciphertexts using \eqref{4a} to obtain $\mathcal{E}(\sum_{i=1}^{n} A_{ui} x_i^k + c)$ and $\mathcal{E}(\sum_{i=1}^{n} A_{gi} x_i^k + d)$. Then the SO sends  $\mathcal{E}(\sum_{i=1}^{n} A_{ui} x_i^k + c)$ and $\mathcal{E}(\sum_{i=1}^{n} A_{gi} x_i^k + d)$ to the agents.

\State  Agents receive and decrypt $\mathcal{E}(\sum_{i=1}^{n} A_{ui} x_i^k + c)$ and $\mathcal{E}(\sum_{i=1}^{n} A_{gi} x_i^k + d)$ from the SO with the private key $(\eta,\mu)$ by using \eqref{2}, then converts the integer to real number by using \eqref{12}.

\State  Each agent $i$ updates the primal variable $x_i$ by \eqref{11a} and the dual variable $\lambda$ by \eqref{11b} over plaintext.

\State Agents calculate the error $\epsilon$, $k=k+1$. 

% \State  $k=k+1$.
\EndWhile
% \EndProcedure
\end{algorithmic}
\label{paradigm_1}
\end{paradigm}

\noindent {\bf{Remark 1:}}
The privacy preservation of the SO is fulfilled via random values $\gamma_i^k$ and $\upsilon_i^k$, i.e., the $i$th agent only receives and decrypts $\mathcal{E}(\gamma_i^kc)$ and $\mathcal{E}(\upsilon_i^kd)$ to get two random numbers $\gamma_i^kc$ and $\upsilon_i^kd$ and therefore acquires no knowledge about the global coefficients $c$ or $d$.
\hfill $\square$

\noindent {\bf{Remark 2:}}
Paradigm \ref{paradigm_1} also allows each agent to generate a separate cryptosystem, i.e., each agent owns a set of public key $(p_i,g_i)$ and private key $(\eta_i,\mu_i)$ \cite{lu2018privacy, hadjicostis2020privacy}. In this case, each agent publicizes $(p_i,g_i)$ and keeps $(\eta_i,\mu_i)$ private. However, such pattern arises significant computation cost and communication cost. Comparing with Paradigm \ref{paradigm_1} where all agents adopt a universal cryptosystem, this  pattern requires each agent to encrypt its message $n$ times with $(p_i,g_i) \forall i \in \mathcal{N}$ and requires the SO to aggregate the collected ciphertexts $n$ times. Therefore, Paradigm \ref{paradigm_1}
achieves privacy preservation with higher computational and communicational efficiency compared to the work in \cite{lu2018privacy, hadjicostis2020privacy}. \hfill $\square$

\section{Privacy Analysis}\label{Section_analysis}

Particularly, we consider three types of adversaries:
\emph{External eavesdroppers} that launch attacks by wiretapping and intercepting exchanged messages between agents and the SO; \emph{SO} which may infer the decision variables of the agents by collecting the received data from the agents;  \emph{Honest-but-curious agents} that follow the paradigm but may collect and observe the received data to infer the private information of other participants. To analyze the privacy, we firstly present the definition of computational indistinguishability. 

\noindent \textbf{Definition 1 \cite{lu2018privacy}:}
Let $\mathcal{X} = \{ X(\kappa)\mid  \kappa \in \mathbb{N} \}$ and 
$\mathcal{Y} = \{Y(\kappa) \mid  \kappa \in \mathbb{N} \}$ be two families of random variables. Let $\mathcal{A}$ be a class of any poly-time algorithms. For every distinguisher $A \in \mathcal{A}$, and for every positive polynomial $p: \mathbb{N} \mapsto \mathbb{R}_{+}$, and every sufficiently large $\kappa$, if it holds that 
\begin{equation}
    |\operatorname{Pr}[A(X(\kappa_{x}))=1]-\operatorname{Pr}[A(Y(\kappa_{y}))=1]|
    < \frac{1}{p(\kappa)},
\end{equation}
where $X(\kappa_{x})$ and $Y(\kappa_{y})$ are random samples drawn from $\mathcal{X}$ and $\mathcal{Y}$, then we say that $\mathcal{X}$ and $\mathcal{Y}$ are computationally indistinguishable, written as 
\begin{equation}
    \mathcal{X}  \stackrel{\text{comp}}{\equiv}\mathcal{Y}.
\end{equation}
\hfill $\blacksquare$

\noindent\textbf{Definition 2 \cite{cramer2015secure}:} Let $\Psi$ be an algorithm for $M$ parties to collectively compute $h_1,\ldots,h_M$, where the $i$th party aims at solving a function $h_i$ that depends on all elements in $x$. We say that $\Psi$ securely computes $h_1,\ldots,h_M$ if there exists a probabilistic poly-time algorithm $S$ such that for the $i$th party with private data $\mathbb{I}_i$, it holds that
\begin{align}
  S(i,\mathbb{I}_i,h_i(\mathbb{I})) \stackrel{\text{comp}}{\equiv} V_i,
\end{align}
where $S(\cdot)$ denotes the overall messages that can be seen after the execution of $S$, $\mathbb{I} \triangleq \{\mathbb{I}_1,\cdots,\mathbb{I}_M\}$, and $V_i$ denotes the data seen by the $i$th party during the execution of $\Psi$. \hfill $\blacksquare$

% Privacy of the agents can be guaranteed if there exists a probabilistic algorithm, defined as a simulator $S(\cdot)$, can produce an output whose distribution is exactly the same as the entire set of data $\mathbb{I}$ seen by the $i$th agent, i.e.,

% \begin{align}
%   S(\mathbb{I}) \stackrel{\text{comp}}{\equiv} \mathbb{I}.
% \end{align} \hfill $\blacksquare$

Definition 1 provides a way to proving
two families of random variables are computationally indistinguishable, and further elaborating on the semantically secure property of a cryptosystem. Suppose an adversary knows two plaintexts $y_1$ and $y_2$ and let an encryption scheme $\bar{\mathcal{E}}(\cdot)$ output two ciphertexts $\bar{\mathcal{E}}(y_1)$ and $\bar{\mathcal{E}}(y_2)$. Then, $\bar{\mathcal{E}}(y_1)$ and $\bar{\mathcal{E}}(y_2)$ are sent to the adversary without telling the adversary which ciphertext corresponds to which plaintext. 
We say that $\bar{\mathcal{E}}(\cdot)$ is semantically secure if for any plaintexts $y_1$ and $y_2$ chosen by the adversary, it holds that $ \bar{\mathcal{E}}(y_1) \stackrel{\text{comp}}{\equiv} \bar{\mathcal{E}}(y_2)$ to the adversary. Definition 2 states that each party could only learn what it should know during the algorithm execution. The security of the proposed privacy-preserving paradigm is presented via Theorem 1 and Theorem 2.

\noindent \textbf{Theorem 1:} Paradigm \ref{paradigm_1} preserves the privacy of the agents against   \emph{external eavesdroppers}, the \emph{SO}, and the \emph{honest-but-curious agents}. \hfill $\blacksquare$

\emph{Proof of Theorem 1:} As shown in \cite{paillier1999public}, the Paillier cryptosystem is semantically secure under decisional composite
residuosity assumption (DCRA), i.e., given a composite $\omega$ and an integer $\phi$, it is computationally hard to decide whether $\phi$ is an $\omega$-residue modulo $\omega^2$. Therefore, the \textit{SO} and \emph{external eavesdroppers} are only able to receive (or wiretap) a sequence of ciphertexts that are computationally indistinguishable, being disabled in extracting useful information about the agents. Then, we consider the privacy preservation of the agents against other \textit{honest-but-curious agents}. Note that the \textit{honest-but-curious agents} could use any received data to infer the privacy of others but do not wiretap the communication links. Let $V_i^k$ denote the set of data that are accessible to the $i$th agent at the $k$th iteration, which can be written as
\begin{align}
  V_i^k &= (A_{ui},A_{gi}, x_i^k,\mathbb{X}_i,\lambda^k,\mathbb{D}, \mathcal{I}_{fi},p,g, \eta,\mu, \gamma_i^kc, \upsilon_i^k d,\nonumber\\
  & \qquad \delta, \mathcal{E}(z_c^k),\mathcal{E}(z_d^k), z_c^k, z_d^k),
\end{align}
where $z_c^k \triangleq \sum_{i=1}^{n} A_{ui} x_i^k + c$, $z_d^k \triangleq \sum_{i=1}^{n} A_{gi} x_i^k + d$, and $\mathcal{I}_{fi}$ denotes the set of coefficients contained in $f_i(x_i)$. Based on Definition 2, we need to construct a poly-time algorithm $S$, such that 
\begin{align}
    S_i^k \stackrel{\text{comp}}{\equiv} V_i^k.
    \label{28s}
\end{align}
Based on \eqref{28s}, agent $i$ then has only to simulate $\mathcal{E}(z_c^k)$ and $\mathcal{E}(z_d^k)$ by generating $\mathcal{E}(z_c^k)^\prime$
and $\mathcal{E}(z_d^k)^\prime$
via $S$ such that 
\begin{subequations}
\begin{align}
   \mathcal{E}(z_c^k) &\stackrel{\text{comp}}{\equiv}  \mathcal{E}(z_c^k)^\prime, \label{29as}\\
   \mathcal{E}(z_d^k) &\stackrel{\text{comp}}{\equiv}  \mathcal{E}(z_d^k)^\prime. \label{29bs}
\end{align}
\end{subequations}
We take \eqref{29as} for example. Eqn. \eqref{29as} holds because 
\begin{align}
    \mathcal{E}(z_c^k) &= \mathcal{E}(\sum_{i=1}^{n} A_{ui} x_i^k +c)\nonumber\\
    &= \prod_{i=1}^{n}  \mathcal{E}(A_{ui} x_i^k +\gamma_i^kc) \nonumber\\
    &= \prod_{i=1}^{n}
    g^{A_{ui} x_i^k +\gamma_i^k c} (r_i^k)^{p} \bmod p^{2}.
\end{align}
From agent $i$'s perspective, each $r_j^k$, $\forall j\in \mathcal{N}$, $j\neq i$ is a random number. Therefore, agent $i$ needs to simulate $\mathcal{E}(z_c^k)^\prime$ by generating ${ (A_{uj} x_j^k + \gamma_j^kc)}^\prime$, $\forall j\in \mathcal{N}$, $j\neq i$, such that 
\begin{align}
    z_c^k = (z_c^k)^\prime.
\end{align}
Then, agent $i$ randomly chooses $(r_j^k)^\prime$ following the same distribution as $r_j^k$ and computes
$\mathcal{E}(A_{ui} x_i^k + \gamma_i^kc)$ and $\mathcal{E}(A_{uj} x_j^k + \gamma_j^kc)$, $\forall j\in \mathcal{N}$, $j\neq i$ as 
\begin{subequations}
\begin{align}
\mathcal{E}(A_{ui} x_i^k + \gamma_i^kc) &= g^{A_{ui} x_i^k +\gamma_i^k c} (r_i^k)^{p} \bmod p^{2}, \label{16as}\\
\mathcal{E}({(A_{uj} x_j^k +\gamma_j^k c)^{\prime}}) &= g^{(A_{uj} x_j^k +\gamma_j^k c)^{\prime}} ((r_j^k)^{\prime})^{p} \bmod p^{2}. \label{16bs}
\end{align}
\end{subequations}
Therefore, agent $i$ is ready to compute
\begin{equation}
    \mathcal{E}(z_c^k)^\prime = 
    \mathcal{E}(A_{ui} x_i^k + \gamma_i^kc) \prod_{j\neq i,j \in \mathcal{N}} \mathcal{E}((A_{uj} x_j^k + \gamma_j^kc)^\prime).
\end{equation}
Further, we can readily have 
\begin{equation}
    \mathcal{D}(\mathcal{E}(z_c^k)^\prime) = z_c^k.
\end{equation}
By following similar lines, \eqref{29bs} can be proved. Therefore, \eqref{28s} is satisfied and Paradigm \ref{paradigm_1} securely computes $\mathcal{E}(z_c^k)$ and $\mathcal{E}(z_d^k)$ between the agents. \hfill $\square$

\noindent \textbf{Theorem 2:} Paradigm \ref{paradigm_1} preserves the privacy of the SO against \emph{external eavesdroppers} and \emph{honest-but-curious agents}. \hfill $\blacksquare$

\emph{Proof of Theorem 2:} During each iteration, the SO sends $\mathcal{E}(\gamma_i^kc)$ and $\mathcal{E}(\upsilon_i^kd)$, which contain the private information $c$ and $d$, to the $i$th agent. Paradigm \ref{paradigm_1}
protects the privacy of the SO against \textit{honest-but-curious agents} by introducing randomized $\gamma_i^k$ and $\upsilon_i^k$ (see also Remark 1). \emph{External eavesdroppers} cannot intercept anything owing to the semantic security. In an extreme case, where all agents are \textit{honest-but-curious agents} and act in collusion, e.g., calculate $\sum_{i=1}^{n}\gamma_i^kc$, the privacy of the SO could be compromised. However, we claim that this case rarely happens in a large-scale network because of the agents' huge population size. \hfill $\square$

% , and it is unavoidable in any privacy-preserving paradigm since each agent must know the true decision variable belonging to itself. \hfill $\square$

\noindent {\bf{Remark 3:}}
Comparing with the private key based algorithms, e.g., SingleMod-based method in \cite{lu2018privacy,huo2020novel}, the proposed public key based paradigm has enhanced security. Lu \emph{et al.} \cite{lu2018privacy} proposed a nonoverlapping partition of the coefficients in the first-order gradients of the Lagrangian to avoid repeated encryption due to the lack of semantic security. However, as the decision variables converge, it is inevitable to encrypt the same decision variables multiple times.  Our preliminary work \cite{huo2020novel} designed a privacy-preserving paradigm 
which naturally avoids the repeated encryption of the coefficients. However, in both \cite{lu2018privacy} and \cite{huo2020novel}, if the SO has access to the update rule of the decision variables, it can use the collected data to estimate the agents’ true decision variables. In contrast, Paradigm \ref{paradigm_1} that is developed based on the Paillier cryptosystem is semantically secure, and the transmitted ciphertexts reveal nothing to the \emph{external eavesdroppers} or the \emph{SO}. \hfill $\square$

\section{Experimental Results}

In this section, we conduct experiments on a highly coupled optimization problem across three agents and one SO. Results on numerical and experimental examples are analyzed to show the efficacy and efficiency of the proposed decentralized privacy-preserving optimization paradigm.

\subsection{Numerical Simulation}
Consider a strongly coupled optimization problem
\begin{equation}
\begin{aligned}
& \underset{x}{\text{min}} & & {\mathcal{F}(x)} \\
& \: \text{s.t.} & & x_{i} \in [0,1], \quad  i=1,2,3,\\
& & &  \sum_{i=1}^{3} A_{gi} x_i + d \leq 0,
\end{aligned}
\label{14}
\end{equation}
where $x_{i} \in \mathbb{R}^2$ and ${\mathcal{F}(x)}$ is the objective function given by 
\begin{align}
    {\mathcal{F}(x)} = &\frac{1}{2}\| \sum_{i=1}^{3} A_{ui} x_i +c \|_2^2
    + \sum_{i=1}^{3}
    (A_{qi}x_i)^{\mathsf{T}}(A_{qi}x_i) \nonumber\\
    &+ \sum_{i=1}^{3}\left(A_{li} x_i + C_{ti} \right).
\end{align}
The coefficients  were chosen as
$c {=} [\begin{smallmatrix} 1\\ 1 \end{smallmatrix}]$,
$A_{u1} {=} [\begin{smallmatrix} -0.2&0\\ 1&-0.5 \end{smallmatrix}]$,
$A_{u2} {=} [\begin{smallmatrix} 0&-2\\ 0&-10 \end{smallmatrix}]$, 
$A_{u3} {=} [\begin{smallmatrix} 1&0\\ 0&1 \end{smallmatrix}]$, $d {=} [\begin{smallmatrix} -1\\ 1 \end{smallmatrix}]$,
$A_{q1} {=} [\begin{smallmatrix} 1&0\\ 1&1 \end{smallmatrix}]$, $A_{q2} {=} [\begin{smallmatrix} 0&1\\ 1&1 \end{smallmatrix}]$,
$A_{q3} {=} [\begin{smallmatrix} -1&1\\ -2&1 \end{smallmatrix}]$,
$A_{l1} {=} [\begin{smallmatrix} 1&1 \end{smallmatrix}]$, $A_{l2} {=} [\begin{smallmatrix} 1&0 \end{smallmatrix}]$,
$A_{l3} {=} [\begin{smallmatrix} 0&1 \end{smallmatrix}]$,
$C_{t1} {=}1$,  $C_{t2} {=}0$, $C_{t3} {=}0.5$, $A_{g1} {=} [\begin{smallmatrix} 0&1\\ 1&-1 \end{smallmatrix}]$, $A_{g2} {=} [\begin{smallmatrix} 0&1\\ -1&1 \end{smallmatrix}]$, and
$A_{g3} {=} [\begin{smallmatrix} 0&0\\ -1&1 \end{smallmatrix}]$. Note that the coefficients were chosen as signed real numbers. The SO possesses all local coefficients $A_{ui}$, $A_{gi}$, $c$, and $d$. $\gamma_i^k$ and $\upsilon_i^k$ are uniformly generated within $[0,1]$ and normalized respectively by the SO. Each agent only learns the local coefficients concerning itself and the individual coefficients set $\mathcal{I}_{fi}=\{A_{qi}, A_{li}, C_{ti}\}$. The coefficients assignment is given by Table \ref{table_numerical_coe_assign}.
% \vspace*{-2mm}
\begin{table}[!htb]
\caption{Coefficients assignment of the numerical example}
\vspace*{-3mm}
\label{table_numerical_coe_assign}
\begin{center}
\begin{tabular}{c|l}
\hline
Participant Name & Coefficients Held \\
\hline
SO & $c,d, A_{ui},A_{gi}, i=1,2,3$ \\
 Agent 1 & $\mathcal{I}_{f1},A_{u1},A_{g1}$    \\
 Agent 2 & $\mathcal{I}_{f2},A_{u2},A_{g2}$   \\
 Agent 3 & $\mathcal{I}_{f3},A_{u3},A_{g3}$   \\
\hline
\end{tabular}
\end{center}
\vspace*{-2mm}
\end{table}

\begin{figure}[!htb]
% \vspace*{-2mm}
    \centering
 \includegraphics[width=0.5\textwidth,trim = 0mm 0mm 0mm 11mm, clip]{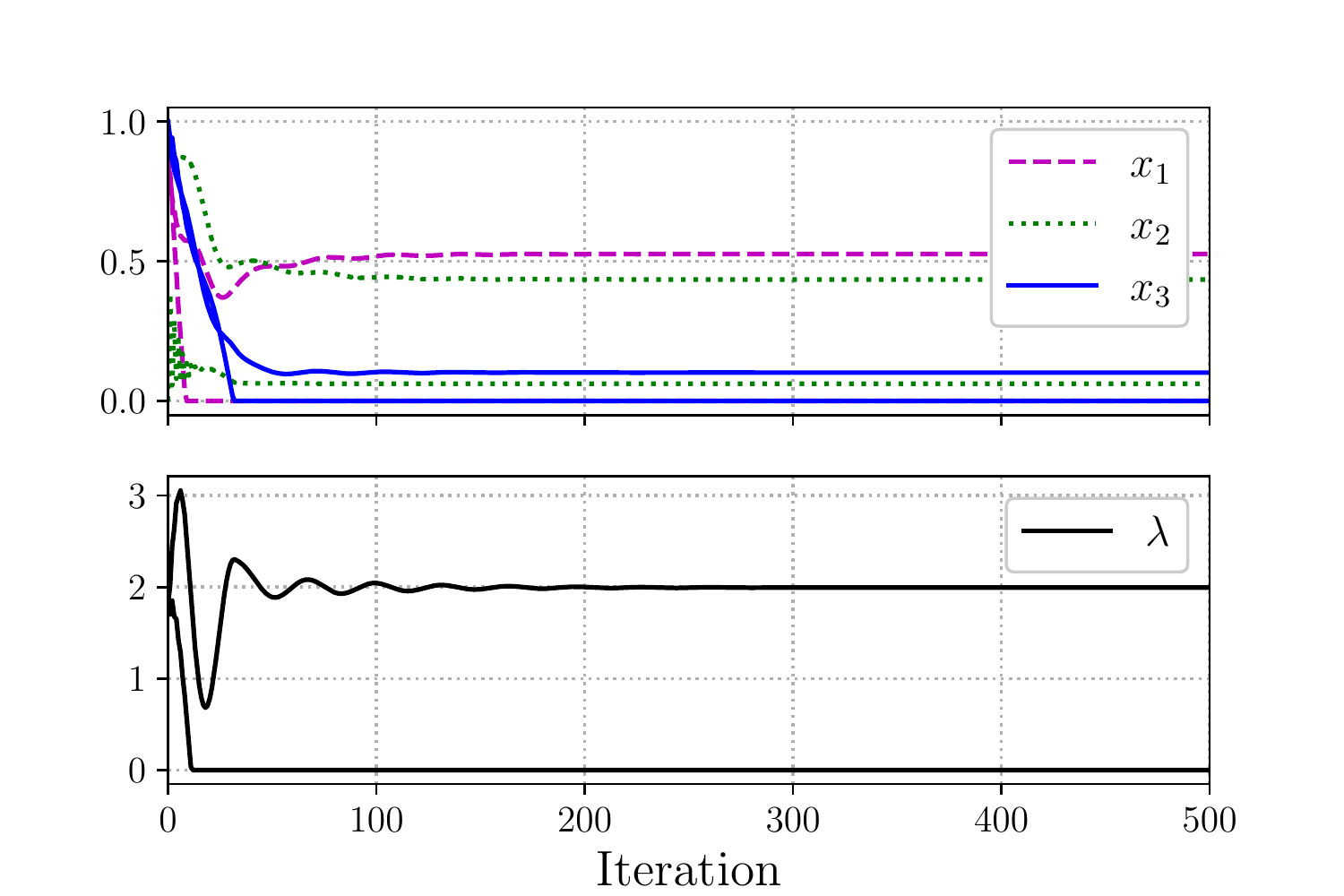}
    \caption{Convergence of the primal and dual variables.}
    \label{numerical_primal_dual_convergence}
\end{figure}

\begin{figure}[!htbp]
% \vspace*{-2mm}
    \centering
 \includegraphics[width=0.5\textwidth,trim = 0mm 0mm 0mm 11mm, clip]{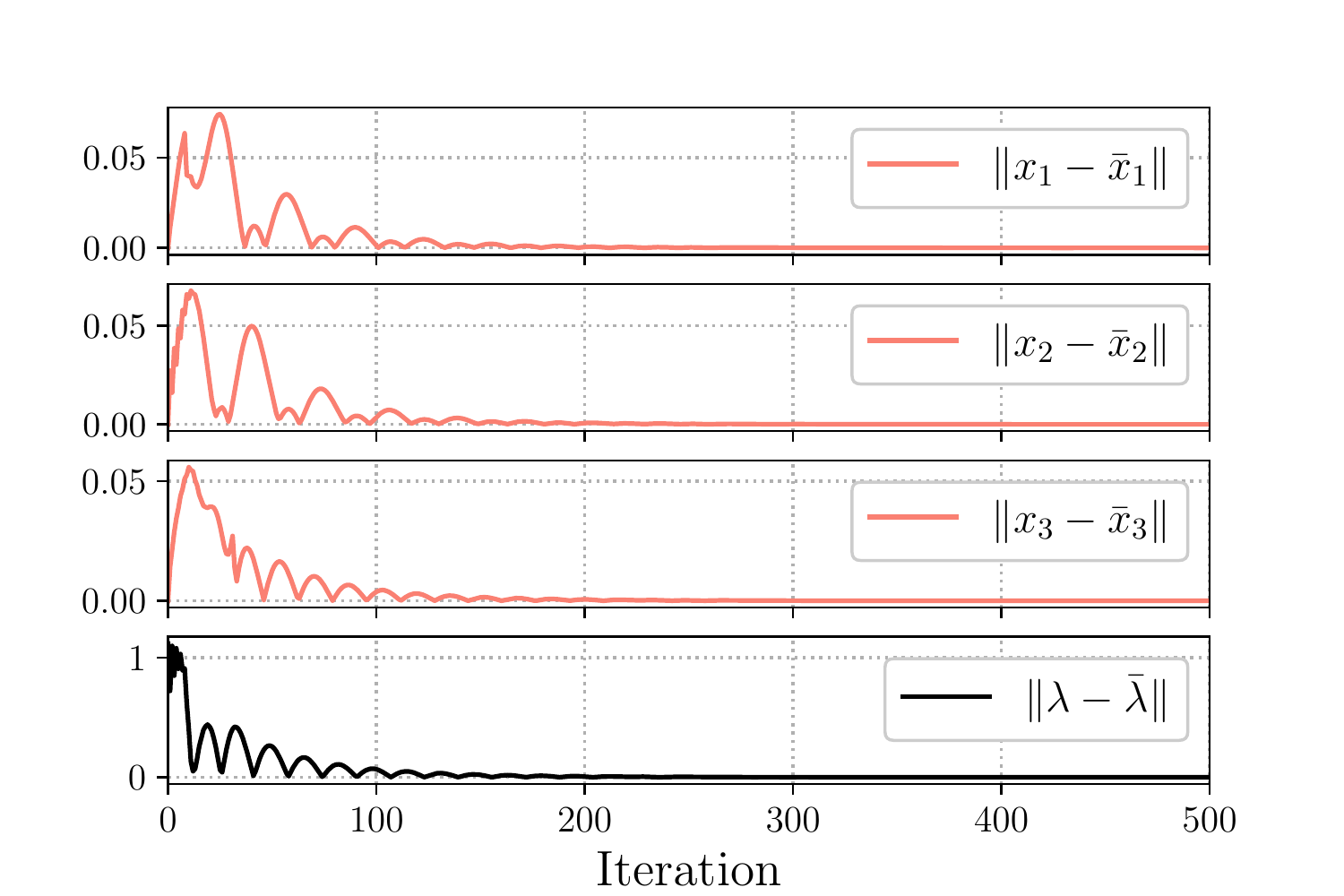}
    \caption{Primal and dual variable gaps between with and without privacy-preservation (precision level $\sigma = 4$).}
    \label{numerical_convergence_error}
\end{figure}

The true optimizers of \eqref{14} are $x_1^* = [0, 0.5258]^{\mathsf{T}}$,
$x_2^* = [0.4347, 0.0621]^{\mathsf{T}}$, and $x_3^* = [0.1016, 0]^{\mathsf{T}}$. We then solve \eqref{14} using Paradigm \ref{paradigm_1}. The primal step sizes were uniformly chosen as $\alpha = 1.6 \times 10^{-2}$ and the dual step size was set to $\beta = 0.8$. Fig. \ref{numerical_primal_dual_convergence}  shows that the primal and dual variables converge in about 200 iterations with $\epsilon_0 = 10^{-3}$.  
Fig. \ref{numerical_convergence_error}  presents the gap of the primal variables and dual variables between the solutions $x_i$ and $ \lambda$ without privacy-preservation and the solutions $\bar{x}_i$ and $ \bar{\lambda}$ of the proposed paradigm. Since only four decimal digits were kept, the proposed approach displayed a slightly slower convergence speed compared to the iterations without privacy-preservation. This can be compensated for by increasing the decimal precision $\sigma$ with extra computational cost.

\subsection{Experimental Platform}
\begin{figure}[!htbp]
% \vspace*{-2mm}
    \centering
    \includegraphics[width=0.5\textwidth,trim = 0mm 0mm 0mm 0mm, clip]{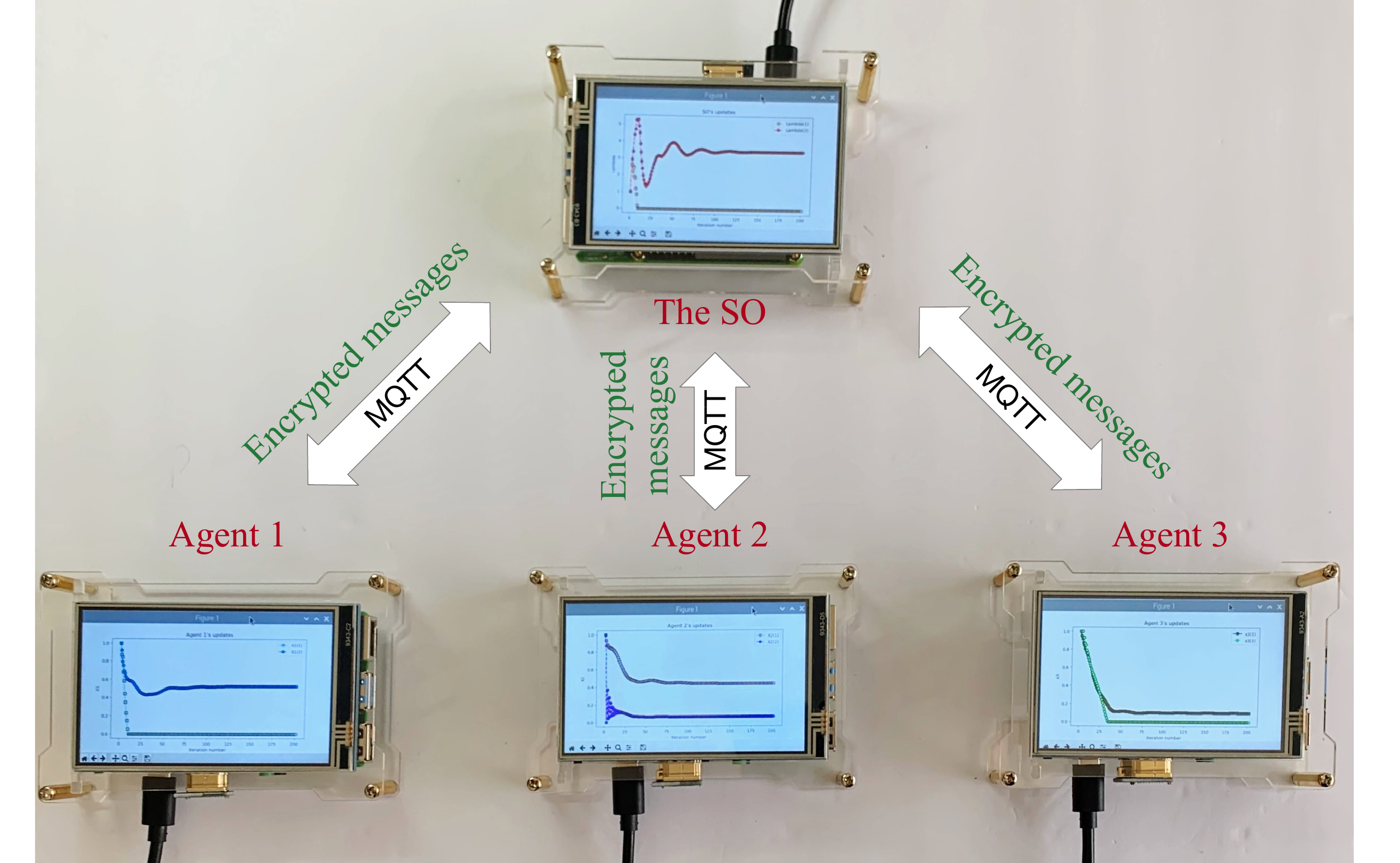}
    % \vspace*{-2mm}
    \caption{Platform structure and communications.}
    \label{physical_platform}
\end{figure}

\begin{figure}[!htbp]
% \vspace*{-2mm}
    \centering
   \includegraphics[width=0.5\textwidth,trim = 8mm 0mm 8mm 8mm, clip]{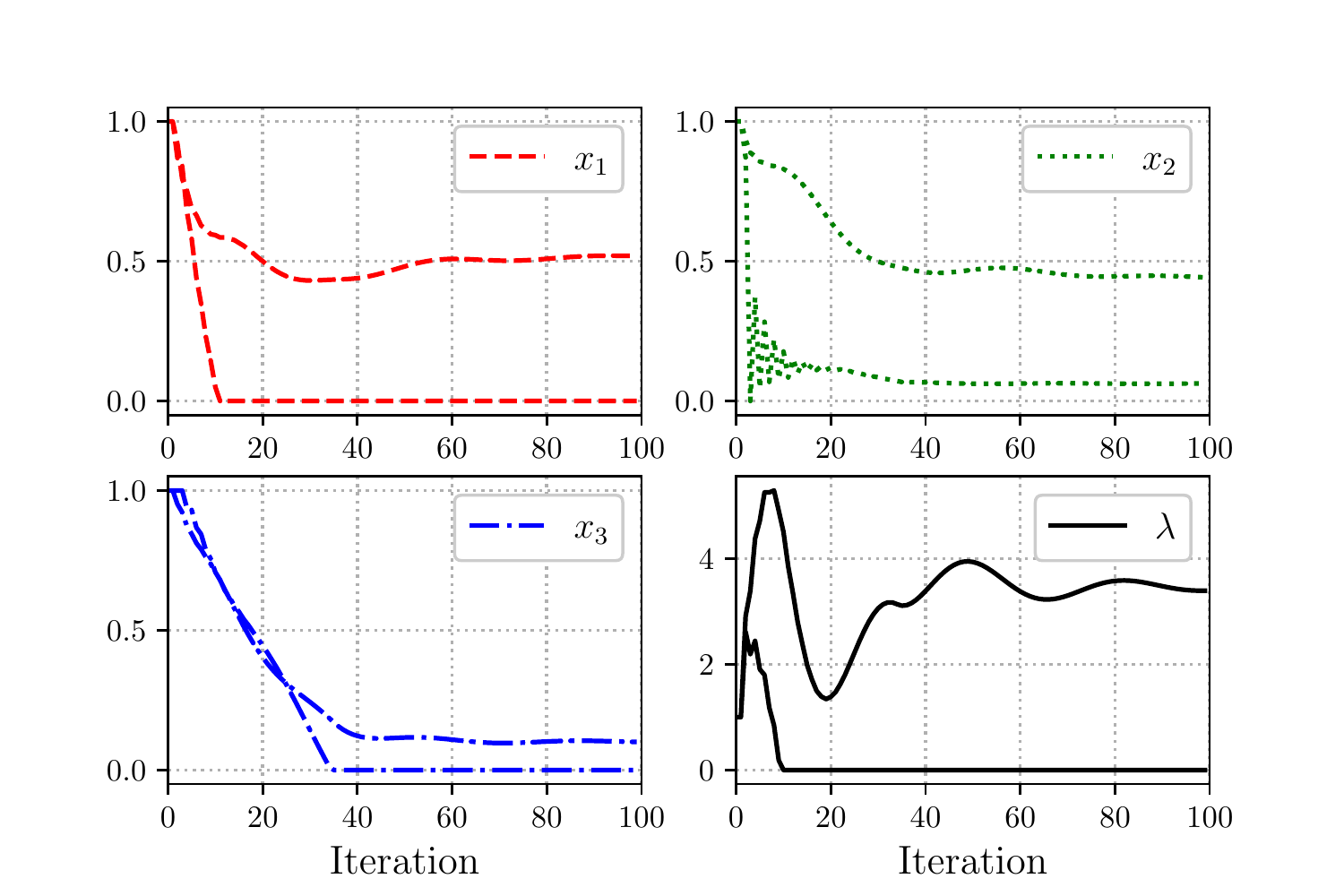}
    % \vspace*{-2mm}
    \caption{Experimental data collected for the convergence of the primal and dual variables.}
    \label{physical_data}
\end{figure}

We also implemented the developed privacy-preserving paradigm on four Raspberry Pi boards to demonstrate the industrial feasibility and value on real-world physical systems. As shown in Fig. \ref{physical_platform}, 
the upper Raspberry Pi represents the SO and the lower three represent the agents. The agents and the SO communicate via Message Queuing Telemetry Transport (MQTT) through Wi-Fi. MQTT can scale to connect with millions of devices, therefore supporting large-scale cooperative optimization problems. The subscribe/publish model is adopted for bi-directional communications. Besides, minimal resources are required for MQTT client, so that the developed paradigm can be implemented on small micro-controllers. The converging processes are dynamically visualized on four LCD screens. Note that the dual variables are updated by agents, while we let the SO plot the dual variables only for clear presentation.  Fig. \ref{physical_data}  presents the real-time data of the convergence for both primal and dual variables (as in Fig. \ref{physical_platform}). We can readily see that the convergence accords with Fig. \ref{numerical_primal_dual_convergence}.

\section{Conclusion}
In this work, we developed a 
decentralized privacy-preserving paradigm to solve a class of strongly coupled multi-agent cooperative optimization problems.
The novel decentralized privacy-preserving paradigm offers scalability, high fidelity, and enhanced security compared to the state of the art.
The private information of both the agents and the SO are protected against the internal and external adversaries, and the privacy preservation is theoretically analyzed and proved. The numerical evaluations proved that the developed paradigm can achieve ideal security without sacrificing fidelity, and the micro-controller based experimental results illustrated the applicability of the proposed paradigm in real-world physical applications.

\addtolength{\textheight}{-12cm}   % This command serves to balance the column lengths

%%%%%%%%%%%%%%%%%%%%%%%%%%%%%%%%%%%%%%%%%%%%%%%%%%%%%%%%%%%%%%%%%%%%%%%%%%%%%%%%

%%%%%%%%%%%%%%%%%%%%%%%%%%%%%%%%%%%%%%%%%%%%%%%%%%%%%%%%%%%%%%%%%%%%%%%%%%%%%%%%

%%%%%%%%%%%%%%%%%%%%%%%%%%%%%%%%%%%%%%%%%%%%%%%%%%%%%%%%%%%%%%%%%%%%%%%%%%%%%%%%
% \section*{Appendix}

% \section*{ACKNOWLEDGMENT}

% The preferred spelling of the word ÒacknowledgmentÓ in America is without an ÒeÓ after the ÒgÓ. Avoid the stilted expression, ÒOne of us (R. B. G.) thanks . . .Ó  Instead, try ÒR. B. G. thanksÓ. Put sponsor acknowledgments in the unnumbered footnote on the first page.

%%%%%%%%%%%%%%%%%%%%%%%%%%%%%%%%%%%%%%%%%%%%%%%%%%%%%%%%%%%%%%%%%%%%%%%%%%%%%%%%

\bibliographystyle{IEEEtran}

\bibliography{bibliography}

\begin{thebibliography}{10}
\providecommand{\url}[1]{#1}
\csname url@samestyle\endcsname
\providecommand{\newblock}{\relax}
\providecommand{\bibinfo}[2]{#2}
\providecommand{\BIBentrySTDinterwordspacing}{\spaceskip=0pt\relax}
\providecommand{\BIBentryALTinterwordstretchfactor}{4}
\providecommand{\BIBentryALTinterwordspacing}{\spaceskip=\fontdimen2\font plus
\BIBentryALTinterwordstretchfactor\fontdimen3\font minus
  \fontdimen4\font\relax}
\providecommand{\BIBforeignlanguage}[2]{{%
\expandafter\ifx\csname l@#1\endcsname\relax
\typeout{** WARNING: IEEEtran.bst: No hyphenation pattern has been}%
\typeout{** loaded for the language `#1'. Using the pattern for}%
\typeout{** the default language instead.}%
\else
\language=\csname l@#1\endcsname
\fi
#2}}
\providecommand{\BIBdecl}{\relax}
\BIBdecl

\bibitem{koshal2011multiuser}
J.~Koshal, A.~Nedi{\'c}, and U.~V. Shanbhag, ``Multiuser optimization:
  Distributed algorithms and error analysis,'' \emph{SIAM Journal on
  Optimization}, vol.~21, no.~3, pp. 1046--1081, 2011.

\bibitem{nimalsiri2019survey}
N.~I. Nimalsiri, C.~P. Mediwaththe, E.~L. Ratnam, M.~Shaw, D.~B. Smith, and
  S.~K. Halgamuge, ``A survey of algorithms for distributed charging control of
  electric vehicles in smart grid,'' \emph{IEEE Transactions on Intelligent
  Transportation Systems}, vol.~21, no.~11, pp. 4497--4515, 2019.

\bibitem{dorronsoro2011improving}
B.~Dorronsoro and P.~Bouvry, ``Improving classical and decentralized
  differential evolution with new mutation operator and population
  topologies,'' \emph{IEEE Transactions on Evolutionary Computation}, vol.~15,
  no.~1, pp. 67--98, 2011.

\bibitem{morstyn2018designing}
T.~Morstyn, A.~Teytelboym, and M.~D. McCulloch, ``Designing decentralized
  markets for distribution system flexibility,'' \emph{IEEE Transactions on
  Power Systems}, vol.~34, no.~3, pp. 2128--2139, 2018.

\bibitem{lu2018privacy}
Y.~Lu and M.~Zhu, ``Privacy preserving distributed optimization using
  homomorphic encryption,'' \emph{Automatica}, vol.~96, pp. 314--325, 2018.

\bibitem{han2016differentially}
S.~Han, U.~Topcu, and G.~J. Pappas, ``Differentially private distributed
  constrained optimization,'' \emph{IEEE Transactions on Automatic Control},
  vol.~62, no.~1, pp. 50--64, 2016.

\bibitem{zhang2016dynamic}
T.~Zhang and Q.~Zhu, ``Dynamic differential privacy for {ADMM}-based
  distributed classification learning,'' \emph{IEEE Transactions on Information
  Forensics and Security}, vol.~12, no.~1, pp. 172--187, 2016.

\bibitem{zhang2018enabling}
C.~Zhang and Y.~Wang, ``Enabling privacy-preservation in decentralized
  optimization,'' \emph{IEEE Transactions on Control of Network Systems},
  vol.~6, no.~2, pp. 679--689, 2018.

\bibitem{wang2019privacy}
Y.~Wang, ``Privacy-preserving average consensus via state decomposition,''
  \emph{IEEE Transactions on Automatic Control}, vol.~64, no.~11, pp.
  4711--4716, 2019.

\bibitem{johnson2019security}
L.~Johnson, \emph{Security Controls Evaluation, Testing, and Assessment
  Handbook}.\hskip 1em plus 0.5em minus 0.4em\relax Academic Press, 2019.

\bibitem{hadjicostis2020privacy}
C.~N. Hadjicostis and A.~D. Dom{\'\i}nguez-Garc{\'\i}a, ``Privacy-preserving
  distributed averaging via homomorphically encrypted ratio consensus,''
  \emph{IEEE Transactions on Automatic Control}, vol.~65, no.~9, pp.
  3887--3894, 2020.

\bibitem{liu2017decentralized}
M.~Liu, P.~K. Phanivong, Y.~Shi, and D.~S. Callaway, ``Decentralized charging
  control of electric vehicles in residential distribution networks,''
  \emph{IEEE Transactions on Control Systems Technology}, vol.~27, no.~1, pp.
  266--281, 2019.

\bibitem{xu2017distributed}
Y.~Xu, T.~Han, K.~Cai, Z.~Lin, G.~Yan, and M.~Fu, ``A distributed algorithm for
  resource allocation over dynamic digraphs,'' \emph{IEEE Transactions on
  Signal Processing}, vol.~65, no.~10, pp. 2600--2612, 2017.

\bibitem{cramer2015secure}
R.~Cramer, I.~B. Damg{\aa}rd, and J.~B. Nielsen, \emph{Secure Multiparty
  Computation}.\hskip 1em plus 0.5em minus 0.4em\relax Cambridge University
  Press, 2015.

\bibitem{paillier1999public}
P.~Paillier, ``Public-key cryptosystems based on composite degree residuosity
  classes,'' in \emph{Proceedings of the International Conference on the Theory
  and Applications of Cryptographic Techniques}, Prague, Czech Republic, May.
  2-6 1999, pp. 223--238.

\bibitem{huo2020novel}
X.~Huo and M.~Liu, ``A novel cryptography-based privacy-preserving
  decentralized optimization paradigm,'' \emph{arXiv:2012.09285 [math.OC]}, pp.
  1--6, 2020.

\end{thebibliography}

% \begin{thebibliography}{99}

% \end{thebibliography}

\end{document}